\documentclass[11pt]{amsart}
\newtheorem{theorem}{Theorem}[section]
\newtheorem{lemma}[theorem]{Lemma}
\newtheorem{cor}[theorem]{Corollary}

\newtheorem*{theoremA}{Theorem A}
\newtheorem*{theoremB}{Theorem B}

\theoremstyle{definition}
\newtheorem{definition}[theorem]{Definition}

\theoremstyle{remark}
\newtheorem{remark}[theorem]{Remark}

\numberwithin{equation}{section}

\begin{document}

\title[Curvature estimates for  submanifolds in  warped products]{Curvature estimates for  submanifolds\\   in  warped products}

\author{L. J. Al\'\i as}
\address{Departamento de Matem\'{a}ticas, \hfill\break
Universidad de Murcia, Campus de Espinardo\hfill\break  30100 Espinardo,
Murcia,   Spain}
\email{ljalias@um.es}

\author{G. P. Bessa, J. F. Montenegro}
\address{Department of Mathematics\hfill\break Universidade Federal do Cear\'a--UFC\hfill\break
Campus do Pici, 60455-760 Fortaleza--CE,
Brazil}
\email{bessa@mat.ufc.br, fabio@mat.ufc.br}

\author{P. Piccione }
\address{Departamento de Matem\'{a}tica, \hfill\break
Instituto de Matem\'atica e Estat\'\i stica\hfill\break
Universidade de S\~ao Paulo\hfill\break Rua do Mat\~ao 1010, Cidade Universit\'aria
\hfill\break S\~ao Paulo, SP 05508-090, Brazil}
\email{piccione.p@gmail.com}
\thanks{This research is a result of the activity developed  within framework of the Programme in Support of Excellence Groups of the Regi\'{o}n de Murcia, Spain by Fundaci\'{o}n S\'{e}neca,  Regional Agency for Science and Technology (Regional Plan for Science and Technology 2007-2010). The first author was partially supported by MEC project PCI2006-A7-0532, MICINN project MTM2009-10418, and Fundaci\'{o}n S\'{e}neca project 04540/GERM/06, Spain.
 The second author was partially supported  by  CNPq-Brazil and  ICTP Associate Schemes.
The third author was partially supported by a CNPq-Brazil. The fourth author is partially supported by Fapesp, grant 2007/03192-7, and by CNPq-Brazil.}


\date{September 2010}

\dedicatory{Dedicated to Keti Tenenblat on occasion of her 65th anniversary, with admiration.}


\begin{abstract}
We give estimates on the intrinsic and the extrinsic curvature of manifolds that are isometrically
immersed as cylindrically bounded submanifolds of warped products. We also address  extensions of the results
in the case of submanifolds of the total space of a Riemannian submersion.
\end{abstract}

\maketitle
\section{Introduction}An important problem in submanifold theory is the isometric immersion problem, this is for  given  two complete  Riemannian manifolds $(M, g_{M})$ and $(N, g_{N})$,  ${\rm dim}(M)<{\rm dim}( N)$ whether   there exists an isometric immersion $\varphi \colon M\hookrightarrow  N$. When $N$ is the Euclidean space, the  Nash Isometric Embedding Theorem answers affirmatively this question, provided the codimension ${\rm dim}(N)-{\rm dim}( M)$ is sufficiently high, see \cite{nash}. For small codimension, here in this paper meaning that ${\rm dim}(N)-{\rm dim}(M)\leq {\rm dim }(M)-1 $,  the answer, in general, depends on the geometries of $M$ and $N$. For instance, on the sectional curvatures of $M$ and $N$, e.g., a classical result of C. Tompkins \cite{tompkins} extended in  a series of papers, by   Chern and Kuiper \cite{chern-kuiper}, Moore, \cite{moore}, O'Neill \cite{oneil},  Otsuki \cite{otsuki} and Stiel \cite{stiel}  can be summarized as follows.
 \begin{theorem} \label{thmTompkins} Let  $\varphi\colon M^{m}\hookrightarrow N^{n}$, $n\leq 2m-1$, be an isometric immersion of a compact  Riemannian $m$-manifold $M$ into  a Cartan-Hadamard $n$-manifold $N$. Then   the sectional curvatures of $M$ and $N$ satisfies \begin{equation}\sup_{M} K_{M}> \inf_{N} K_{N}. \label{eqTompkins}\end{equation}
 \end{theorem} Theorem \ref{thmTompkins} was  extended,  in a seminal paper \cite{JoKo80}, by L. Jorge and D. Koutrofiotis  to bounded, complete  submanifolds with scalar curvature bounded below  where they introduced the Omori-Yau maximum principle to tackle this immersion problem. Due to a much better version of the Omori-Yau maximum principle, Pigola-Rigoli-Setti \cite{PiRiSe05}, extended Jorge-Kourtofiotis result to   to bounded, complete  submanifolds with scalar curvature satisfying \begin{equation} s_{M}(x)\geq - \rho^{2}_{M}(x)\cdot\prod_{j=1}^{k}\big[\log^{(j)}(\rho_{M}(x))\big]^{2},\,\, \rho_{M}(x)\gg 1 \label{eqScalar}\end{equation} where $\rho_{M}$ is the distance function on $M$ to a fixed point and $\log^{(j)}$ is the $j$-th iterate of the logarithm.
 \begin{theorem}[Jorge-Koutrofiotis-Pigola-Rigoli-Setti]\label{thmJorge-koutrofiotis}Let $\varphi \colon M \hookrightarrow N$ be an isometric immersion of a complete $m$-dimensional Riemannian manifold $M$ into an $n$-dimensional Riemannian manifold $N$, $n\leq 2m-1$.  Let $B_{N}(r)$ be a geodesic ball of $N$ centered at $p=\varphi(q)$ with radius $r<\min\{{\rm inj}_N(p), \pi/2\sqrt{b}\}$, where  the radial sectional
curvatures $K_N^{\mathrm{rad}}$ along the radial geodesics issuing from
$p$ are bounded
as $K_N^{\mathrm{rad}}\leq b$ in $B_N(r)$ and where  $\pi/2\sqrt{b}$ is replaced  by $+\infty$ if $b\leq 0$. Suppose that the scalar curvature of  $M$ satisfies \eqref{eqScalar} and that  $
\varphi (M)\subset B_N(r)$. Then
\begin{equation}\label{eqJorge-koutrofiotis}\sup_{M}K_{M}\geq C_{b}^{2}(r)+\inf_{B_{N}(r)}K_{N}.
\end{equation}Where    $$
C_b(t)
=\left\{\begin{array}{lll}
\sqrt{b}\cot(\sqrt{b}\, t) & \mathrm{if} & b >0,\\
1/t & \mathrm{if} & b =0,\\
\sqrt{-b}\coth(\sqrt{-b}\, t) & \mathrm{if}  & b <0.
\end{array}\right.
$$

\end{theorem}
Very recently,  Theorem \ref{thmJorge-koutrofiotis} was extended to the class of cylindrically bounded  submanifolds  by Alias-Bessa-Montenegro \cite{AlBeMo09}. \begin{theorem}[Alias-Bessa-Montenegro]\label{thmalias-bessa-mont} Let $\varphi\colon\! M^m\hookrightarrow N^{n}\!\times\mathbb{R}^{\ell}$  be  an isometric immersion of a complete Riemannian $m$-manifold $M$ into the product $N^{n}\!\times\mathbb{R}^{\ell}$, $ n+2\ell\leq 2m-1$,  where $N^{n}$ is a  Riemannian $n$-manifold. Let $B_{N}(r)$ be a geodesic ball of $N$ centered at $p=\pi_{1}(\varphi (q))$\footnote{$\pi_{1}\colon N\times \mathbb{R}^{\ell}\to N$ is the projection on the first factor.}    with radius $r<\min\{{\rm inj}_N(p), \pi/2\sqrt{b}\}$, where  the radial sectional
curvatures $K_N^{\mathrm{rad}}$ along the radial geodesics issuing from
$p$ are bounded
as $K_N^{\mathrm{rad}}\leq b$ in $B_N(r)$ and  $\pi/2\sqrt{b}$ is replaced  by $+\infty$ if $b\leq 0$.
 Suppose that
$
\varphi (M)\subset B_N(r)\times \mathbb{R}^\ell
$
 and one of these two conditions below holds,
\begin{enumerate}
\item[i.] $\varphi$ is proper and $\sup_{B_{N}(r)\times B_{\mathbb{R}^{\ell}}(t)} \Vert \alpha \Vert \leq G(t)$, where $G\colon [0, \infty)\to (0, \infty)$ with $1/G \not \in L^{1}(0, \infty)$.
    \item[]
\item[ii.] The scalar curvature of $M$ satisfies $(\ref{eqScalar})$.
\end{enumerate}then
\begin{equation}\sup_{M}K_{M}\geq C_{b}^{2}(r)+\inf_{B_{N}(r)}K_{N}.\label{eq-alias-bessa-montenegro}
\end{equation}

\end{theorem}The purpose of this paper is to show that these results above can be  extended naturally to isometric immersion into warped product manifolds  $\mathcal M=\mathcal X\times\mathcal  V$ endowed with the metric $g^\mathcal M=g^\mathcal X+\psi^2g^\mathcal V$, where $(\mathcal X,g^\mathcal X)$ and $(\mathcal V,g^\mathcal V)$ are  Riemannian  manifolds
and  $\psi\colon\! \mathcal X\!\to\!\mathbb R^+$ is a smooth positive
function on $\mathcal X$. Set $n_{\mathcal X}\!=\!\mathrm{dim}(\mathcal X)$, $n_\mathcal V\!=\!\mathrm{dim}(\mathcal V)$
and $n_\mathcal M\!=\!n_\mathcal X\!+\!n_\mathcal V\!=\!\mathrm{dim}(\mathcal M)$.
We will assume in our main result a certain bound on the radial curvature along geodesics issuing from a point
$x_0$ in the base manifold $\mathcal X$, see \eqref{eq:assradialcurvature}.

Our first result is the following (see Theorem~\ref{thm:mainsectcurvestimate}).
\begin{theoremA}\label{thm:mainsectcurvestimate-introd}Let $(M,g^M)$ be a complete Riemannian $n_{M}$-manifold for which the
weak Omori--Yau principle for the Hessian\footnote{See Section \ref{sec:OmoriYau} for the details  of the weak Omori-Yau maximum principles for the Hessian and for the Laplacian} holds and
let $\varphi:M\to\mathcal M$ be an isometric immersion.
Assume that the following hypotheses are satisfied
\begin{enumerate}
\item\label{itm:sectcurv2-introduction} $\pi_\mathcal X\big(\varphi(M)\big)\subset B_\mathcal X(r)$, a geodesic ball in $\mathcal X$  with center at some $x_0\in\mathcal X$ and $r\in\left(0,
\mathrm{inj}_\mathcal X(x_0)\right)$.
\smallskip

\item Assumption \eqref{eq:assradialcurvature} holds.
\smallskip

\item $2n_M\ge 2n_\mathcal V+n_\mathcal X+1$.
\end{enumerate}
Then,
\begin{equation}\label{eq:sectcurvestimate-introduction}
\sup_M K_{M}\ge C_b(r)^2+\inf_{B_\mathcal X(x_0;r)}K_\mathcal X.
\end{equation}
\end{theoremA}
Our second main result gives an estimate on the mean curvature of cylindrically bounded submanifolds
of warped product (see Theorem~\ref{thm:meancurvimmwarpedpr}).
\begin{theoremB}
Let $(M,g^M)$ be a complete  Riemannian manifold for which the weak Omori--Yau principle for the Laplacian holds, and
let $\varphi:M\to\mathcal M$ be an isometric immersion.
Assume that the following hypotheses are satisfied.
\smallskip
\begin{enumerate}
\item $\pi_\mathcal X\big(\varphi(M)\big)\subset B_\mathcal X(x_0;r)$ a geodesic ball in $\mathcal X$  with center at some $x_0\in\mathcal X$ and $r\in\left(0,
\mathrm{inj}_\mathcal X(x_0)\right)$.
\smallskip

\item Assumption~\eqref{eq:assradialcurvature} holds.
\end{enumerate}
Then, denoting by $\vec H^\varphi$ the mean curvature vector of $\varphi$, one has the following estimate
on the supremum of $\vert\vec H^\varphi\big\vert_\mathcal M$.
\begin{equation}
\sup_M\big\vert\vec H^\varphi\big\vert_\mathcal M\ge(n_M-n_\mathcal V)\, C_b(r)-{n_\mathcal V}\Psi_0,
\end{equation}
where
\begin{equation}
\Psi_0=\sup_{\mathrm{dist}(x,x_0)\le r}\Big\vert\frac{\mathrm{grad}^\mathcal X\psi(x)}{\psi(x)}\Big\vert_\mathcal X.
\end{equation}
\end{theoremB}
In view of the above results, it is an interesting question to establish
geometric conditions for the validity of the
weak Omori-Yau maximum principles in submanifolds of warped products.
We study this question in Section~\ref{sec:OmoriYau}, see Theorem~\ref{thm:OYimmersion}.

We also observe that the results of Theorems A and B can be generalized to the more general situation
of cylindrically bounded isometric immersions into the total space of a Riemannian submersion.
These generalizations are discussed in Section~\ref{sec:immersionssubmersions}, see Theorem~\ref{thm:mainsectcurvestimatesub}
and Theorem~\ref{thm:meancurvimmsub}. In this situation, the curvature estimates are given in terms of the norms
of the characteristic tensors of the submersion.

\section{Generalities on warped products}
Let $(\mathcal X,g^\mathcal X)$ and $(\mathcal V,g^\mathcal V)$ be  Riemannian  manifolds
and let $\psi\colon\! \mathcal X\!\to\!\mathbb R^+$ be a smooth positive
function on $\mathcal X$. Set $n_{\mathcal X}\!=\!\mathrm{dim}(\mathcal X)$, $n_\mathcal V\!=\!\mathrm{dim}(\mathcal V)$
and $n_\mathcal M\!=\!n_\mathcal X\!+\!n_\mathcal V\!=\!\mathrm{dim}(\mathcal M)$.
The product manifold $\mathcal M=\mathcal X\times\mathcal  V$ endowed with the metric $g^\mathcal M=g^\mathcal X+\psi^2g^\mathcal V$ is
the \emph{warped product of $\mathcal X$ and $\mathcal V$}, with \emph{warping function $\psi$}. This is also denoted
with the symbol $\mathcal M=\mathcal X\times_\psi\mathcal V$. The projection $\pi_\mathcal X:\mathcal M\to\mathcal X$
is a Riemannian submersion, while the projection $\pi_\mathcal V:\mathcal M\to\mathcal V$ is not (unless
$\psi\equiv 1$). In fact, warped products are special cases of Riemannian submersions, characterized by
the property of having integrable horizontal distribution with totally geodesic leaves and with totally
umbilical fibers, see Section~\ref{sec:immersionssubmersions}. We will borrow some terminology from Riemannian submersions, and we will call
$\mathcal X$ the \emph{base} and $\mathcal V$ the \emph{fiber} of $\mathcal M$. Moreover, vectors
that are in the kernel of $\mathrm d\pi_\mathcal X$ are called \emph{vertical},
while vectors in the kernel of $\mathrm d\pi_\mathcal V$ are called \emph{horizontal}.

Vector fields $X\in\mathfrak X(\mathcal X)$ will be identified with vector fields
in $\mathcal M$ that ``do not depend on the second variable'', i.e., $X(x,v)=X(x)$ for all $v\in\mathcal V$.
This type of horizontal vector fields will be called \emph{h-basic}.
Similarly, vector fields $V\in\mathfrak X(\mathcal V)$ will be identified with vector fields in $\mathcal M$
that do not depend on the first variable; they will be called \emph{v-basic}.
Note that
\begin{equation}\label{thm:piXisom}
g^\mathcal M(X,Y)=g^\mathcal X(X,Y)
\end{equation}
for every h-basic $X$ and $Y$, while
\begin{equation}\label{thm:piVnotisom}
g^\mathcal M(V,W)=\psi^2g^\mathcal V(V,W)
\end{equation}
for every v-basic $V$ and $W$. We will use consistently the notation $X$, $Y$, $Z$ for
h-basic fields, and $U$, $V$, $W$ for v-basic fields on $\mathcal M$.
Observe that the Lie bracket $[X,Y]$ of h-basic vector fields is h-basic, the Lie bracket
$[V,W]$ of v-basic fields is v-basic, while the Lie bracket $[X,V]$ of an h-basic and a v-basic field
is zero.
\subsection{Riemannian differential operators}
The symbols
$\nabla^\mathcal X$ and $\nabla^\mathcal V$ will denote the Levi--Civita connections
of $(\mathcal X,g^\mathcal X)$ and of $(\mathcal V,g^\mathcal V)$ respectively.
Due to well known invariance properties, the Levi--Civita connection
$\nabla^\mathcal M$ of $(\mathcal M,g^\mathcal M)$ will be uniquely determined by the values
of $\nabla^\mathcal M_AB$, where $A$ and $B$ are basic. Since $[X,V]=0$, then
\begin{equation}\label{eq:LC0}
\nabla^\mathcal M_XV=\nabla^\mathcal M_VX,
\end{equation}
for all h-basic $X$ and v-basic $V$.
The following formulas are easily computed
\begin{equation}\label{eq:LC1}
\nabla^\mathcal M_XY=\nabla^\mathcal X_XY,
\end{equation}
for every h-basic fields $X$ and $Y$. It follows in particular that $\mathcal X\times\{v\}$ is a totally
geodesic submanifold of $\mathcal M$ for all $v\in\mathcal V$; moreover, the curvature tensor
$R^\mathcal M$ of horizontal vector is given by
\[R^\mathcal M(X,Y)Z=R^\mathcal X(X,Y)Z.\]
Thus, the sectional curvature $K_\mathcal M(X,Y)$ of the plane in $T\mathcal M$
spanned
by horizontal vectors $X$ and $Y$ coincides with the sectional curvature in the base manifold
\begin{equation}\label{eq:equalitysectcurvature}
K_\mathcal M(X,Y)=K_\mathcal X(X,Y).
\end{equation}
As to the covariant derivative of mixed terms
\begin{eqnarray}\label{eq:LC2}
&&\nabla^\mathcal M_VX\stackrel{\text{by \eqref{eq:LC0}}}{=}\nabla^\mathcal M_XV=\frac{X(\psi)}\psi\,V,
\\
\label{eq:LC3}
&&\nabla^\mathcal M_VW=\nabla^\mathcal V_VW-g^\mathcal M(V,W)\,\frac{\mathrm{grad}^\mathcal X\psi}\psi,
\end{eqnarray}
for all h-basic $X$ and all v-basic $V$ and $W$. The second fundamental form of the fibers
$\{x\}\times\mathcal V$ is
\begin{equation}\label{eq:secfundformfibers}
\mathcal S^\mathcal V(V,W)=-g^\mathcal M(V,W)\,\frac{\mathrm{grad}^\mathcal X\psi}\psi.
\end{equation}
Therefore, critical points of $\psi$ correspond to totally geodesic fibers. By taking trace
in \eqref{eq:secfundformfibers}, we get the following expression for the mean curvature vector
$\vec H$ of the fibers
\begin{equation}\label{eq:meancurvaturefibers}
\vec H=-n_{\mathcal V}\frac{\mathrm{grad}^\mathcal X\psi}\psi.
\end{equation}
For an h-basic vector field $X$ we have
\begin{equation}\label{eq:divhbasic}
\mathrm{div}^\mathcal M(X)=
\mathrm{div}^\mathcal X(X)+n_{\mathcal V}\,\displaystyle\frac{X(\psi)}{\psi},
\end{equation}
while for a v-basic field $V$
\begin{equation}\label{eq:divvbasic}
\mathrm{div}^\mathcal M(V)=\mathrm{div}^\mathcal V(V).
\end{equation}
Let $F:\mathcal X\to\mathbb R$ be a smooth function and denote by $F^\mathrm h=F\circ\pi_{\mathcal X}:\mathcal M\to\mathbb R$ the lifting of $F$ to $\mathcal M$. It is easily seen that the gradient $\mathrm{grad}^\mathcal MF^\mathrm h$ is
horizontal. The gradient of $F^\mathrm h$ is the h-basic field
\begin{equation}\label{eq:gradh}
\mathrm{grad}^\mathcal MF^\mathrm h=\mathrm{grad}^\mathcal XF.
\end{equation}
Similarly, if $G:\mathcal V\to\mathbb R$ is a smooth function, and $G^\mathrm v=G\circ\pi_{\mathcal V}$ is
its lifting to $\mathcal M$, then the gradient $\mathrm{grad}^\mathcal MG^\mathrm v$ is
vertical, but not v-basic
\begin{equation}\label{eq:gradv}
\mathrm{grad}^\mathcal MG^\mathrm v=\frac1{\psi^2}\,\mathrm{grad}^\mathcal VG.
\end{equation}

The Laplacian $\Delta^\mathcal M$
of the functions $F^\mathrm h$ and $G^\mathrm v$ is given by
\begin{eqnarray}\label{eq:laplacianh}
&\Delta^\mathcal MF^\mathrm h=\Delta^\mathcal XF
+n_\mathcal V\cdot g^\mathcal X\big(\mathrm{grad}^\mathcal XF,\displaystyle\frac{\mathrm{grad}^\mathcal X\psi}{\psi}\big).
\\
\label{eq:laplacianv}
&\Delta^\mathcal MG^\mathrm v=\displaystyle\frac{1}{\psi^2}\,\Delta^\mathcal VG.
\end{eqnarray}
As to the Hessian of the functions $F^\mathrm h$ and $G^\mathrm v$, the following formulas can be computed
easily:
\begin{eqnarray}\label{eq:hessianFhhh}
&&\mathrm{Hess}^\mathcal MF^\mathrm h(X,X)=\mathrm{Hess}^\mathcal XF(X,X),
\\
\label{eq:hessianFhvv}
&&\mathrm{Hess}^\mathcal MF^\mathrm h(V,V)=
g^\mathcal X\big(\mathrm{grad}^\mathcal XF,\frac{\mathrm{grad}^\mathcal X\psi}{\psi}\big)\cdot g^\mathcal M(V,V),
\\
\label{eq:hessianFhhv}
&&\mathrm{Hess}^\mathcal MF^\mathrm h(X,V)=0,
\\\label{eq:hessianGvhh}
&&\mathrm{Hess}^\mathcal MG^\mathrm v(X,X)=0,
\\
\label{eq:hessianGvvv}
&&\mathrm{Hess}^\mathcal MG^\mathrm v(V,V)=\mathrm{Hess}^\mathcal VG(V,V),
\\
\label{eq:hessianGvhv}
&&\mathrm{Hess}^\mathcal MG^\mathrm v(X,V)=-\frac{X(\psi)}{\psi}\,V(G).
\end{eqnarray}
\subsection{Isometric immersions into warped products}
Let us now consider an immersion $\varphi:M\hookrightarrow\mathcal M$. Assume that $M$ is endowed
with the pull-back metric $g^M=\varphi^*(g^\mathcal M)$. If $L:\mathcal M\to\mathbb R$ is a smooth function,
then setting $f=L\circ\varphi:M\to\mathbb R$, one computes\footnote{%
Obviously, formulas \eqref{eq:gradcomposition} and \eqref{eq:hessiancomposition} hold for an isometric immersion
$\varphi:M\to\mathcal M$ into \emph{any} ambient manifold $\mathcal M$, not just in warped products.}
\begin{equation}\label{eq:gradcomposition}
g^M(\mathrm{grad}^Mf,e)=g^\mathcal M\big(\mathrm{grad}^\mathcal ML,\mathrm e\big),
\end{equation}
for all $e\in TM$ and
\begin{equation}\label{eq:hessiancomposition}
\mathrm{Hess}^Mf(e,e)=\mathrm{Hess}^\mathcal ML\big(\mathrm e,\mathrm e\big)+g^\mathcal M\big(
\mathrm{grad}^\mathcal ML,\mathcal S^\varphi(e,e)\big),
\end{equation}
for all $e\in TM$. Here, $\mathcal S^\varphi$ is the second fundamental form of $\varphi$ and  we identified $e$ with $d\varphi(e)$.
In particular, for $L=G^\mathrm v$, in the notation above, using \eqref{eq:gradv}, \eqref{eq:hessianGvhh}, \eqref{eq:hessianGvhv}
and \eqref{eq:hessianGvvv} one has
\begin{multline}\label{eq:hessianimmersionvertical-1}
\mathrm{Hess}^Mf(e,e)=
-2\,g^\mathcal X\big(\frac{\mathrm{grad}^\mathcal X\psi}{\psi},\mathrm e^\mathrm{hor}\big)\,g^\mathcal
V\big(\mathrm{grad}^\mathcal VG,\mathrm e^\mathrm{ver}\big)\\+
\mathrm{Hess}^\mathcal VG\big(\mathrm e^\mathrm{ver},\mathrm e^\mathrm{ver}\big)
+g^\mathcal V\big(\mathrm{grad}^\mathcal VG,\mathcal S^\varphi(e,e)^\mathrm{ver}\big).
\end{multline}

Here, $\mathrm e^\mathrm{hor}$ and $\mathrm e^\mathrm{ver}$ are respectively the horizontal and the
vertical components of $\mathrm e$.
Similarly, for $L=F^\mathrm h$, using \eqref{eq:gradh}, \eqref{eq:hessianFhhh},
\eqref{eq:hessianFhvv} and \eqref{eq:hessianFhhv} one computes
\begin{multline}\label{eq:hessianimmersionhorizontal}
\mathrm{Hess}^Mf(e,e)=
g^\mathcal X\big(\mathrm{grad}^\mathcal XF,\frac{\mathrm{grad}^\mathcal X\psi}\psi\big)\cdot
g^\mathcal M\big(\mathrm e^\mathrm{ver},\mathrm e^\mathrm{ver}\big) \\+ \mathrm{Hess}^\mathcal XF\big(\mathrm e^\mathrm{hor},\mathrm e^\mathrm{hor}\big)
+g^\mathcal X\big(\mathrm{grad}^\mathcal XF,\mathcal S^\varphi(e,e)^\mathrm{hor}\big).
\end{multline}

\begin{section}{On the Omori-Yau Maximum Principle}
\label{sec:OmoriYau}
\begin{definition}[Pigola-Rigoli-Setti]
Let $(M,g^M)$ be a Riemannian manifold. We say that the \emph{Omori--Yau Maximum Principle for the Hessian holds
in $(M,g^M)$}
if for every smooth function $f\colon M\to\mathbb R$ with $\sup\limits_Mf<+\infty$
there exists a sequence $(p_n)_{n\in\mathbb N}$ in $M$ such that
\begin{itemize}
\item[(a)] $\lim\limits_{n\to\infty}f(p_n)=\sup\limits_Mf$,
\smallskip

\item[(b)] $\big\Vert\mathrm{grad}^Mf(p_n)\big\Vert\le\displaystyle\frac{1}{n}$,
\smallskip

\item[(c)]  $\mathrm{Hess}^Mf(p_n)(e,e)\le\displaystyle\frac{1}{n}\,g^M(e,e)$ for
all $e\in T_{p_n}M$,
\end{itemize}
for all $n$. Similarly, the \emph{Omori--Yau Maximum Principle for the Laplacian holds
in $(M,g^M)$}
if the above properties hold, with (c) replaced by the  condition
\begin{itemize}
\item[(c')] $\Delta^Mf(p_n)\le\displaystyle\frac{1}{n}$.
\end{itemize}
We say that the \emph{Weak Omori--Yau Principle for the Hessian} (\emph{Laplacian})
in $(M,g^M)$ if for every smooth function $f\colon M\to\mathbb R$ with $\sup_M f<+\infty$
there exists a sequence $(p_n)_{n\in\mathbb N}$ in $M$ satisfying (a) and (c) ((a) and (c'))
above.
\end{definition}
The following Theorem, due to Pigola, Rigoli and Setti \cite{PiRiSe05}, gives sufficient conditions for  the Omori-Yau Maximum Principle to hold in a Riemannian manifold.
\begin{theorem}[Pigola-Rigoli-Setti]\label{thm:sufcondOY}Let $(M, g^{M})$ be a Riemannian manifold.
Assume that there exist smooth functions \[h\colon \left[0,+\infty\right)\to\left[0,+\infty\right)\quad\text{and}\quad
\gamma\colon M\to\left[0,+\infty\right)\] such that
\begin{enumerate}
\item\label{itm:sufcondOY1} $h(0)>0$ and $h'(t)\ge 0$ for all\, $t\ge0$,
\item[]
\item\label{itm:sufcondOY2} $\limsup\limits_{t\to+\infty}\,\displaystyle t\cdot h\big(\sqrt t\big)/h(t)<+\infty$,
\item[]

\item\label{itm:sufcondOY3} $\displaystyle\int_0^{+\infty}\mathrm dt/\sqrt{\strut h(t)}=+\infty$,
\item[]

\item\label{itm:sufcondOY4} $\gamma$ is proper,
\item[]
\item\label{itm:sufcondOY5} $\big\vert\mathrm{grad}^M\gamma\big\vert\le c\cdot\sqrt\gamma$\, for some $c>0$
outside a compact subset of $M$,
\item[]
\item\label{itm:sufcondOY6} $\mathrm{Hess}^M\gamma\le c'\cdot\sqrt{\strut\gamma\cdot h(\sqrt\gamma)}$\, for some
$c'>0$ outside a compact subset of $M$.
\end{enumerate}
Then the Omori--Yau Maximum Principle for the Hessian holds in $(M,g^M)$.
A totally analogous statement holds in the case of the Omori--Yau principle for the Laplacian, with
assumption \eqref{itm:sufcondOY6} replaced by
\begin{itemize}\item[]
\item[(\ref{itm:sufcondOY6}')] $\triangle^M\gamma\le c'\cdot\sqrt{\strut\gamma\cdot h(\sqrt\gamma)}$ for some
$c'>0$ outside a compact subset of $M$.
\end{itemize}
\end{theorem}

Note that any function $h$ satisfying \eqref{itm:sufcondOY1} and \eqref{itm:sufcondOY2} is unbounded
\[\lim_{t\to+\infty}h(t)=+\infty.\]
\begin{definition}A pair of functions $(h,\gamma)$ satisfying
\eqref{itm:sufcondOY1}---\eqref{itm:sufcondOY6} of Theorem~\ref{thm:sufcondOY} will
be called an \emph{OY-pair} for the Hessian in $(M,g^M)$. Similarly, a pair $(h,\gamma)$ satisfying
\eqref{itm:sufcondOY1}---\eqref{itm:sufcondOY5} and (\ref{itm:sufcondOY6}') is called
an {\em OY-pair} for the Laplacian in $(M,g^M)$.
\end{definition}

\noindent We will now assume that $(M,g^M)$ is isometrically immersed in a warped product $\mathcal X\times_\psi\mathcal V$,
and we want to give conditions that guarantee the validity of the Omori-Yau Maximum Principle on $(M,g^M)$
in terms of the corresponding property of $(\mathcal V,g^\mathcal V)$ and the geometry of the immersion.
\begin{theorem}\label{thm:OYimmersion}
Let $\varphi\colon M\to\mathcal M=\mathcal X\times_\psi\mathcal V$ be an isometric immersion, and let
$(h,\Gamma)$ be an OY-pair for the Hessian in $(\mathcal V,g^\mathcal V)$. Set $\gamma=\Gamma^\mathrm v\circ\varphi\colon M\to\left[0,+\infty\right.)$.
Assume the following hypothesis.
\begin{itemize}
\item[(a)] $\varphi$ is proper,
\item[]
\item[(b)] $\pi_\mathcal X\big(\varphi(M)\big)$ is contained in a compact subset $K$ of $\mathcal X$,
\item[]
\item[(c)] $\big\Vert\mathcal S^\varphi\big\Vert\le\alpha\,\sqrt{\strut h(\sqrt\gamma)}$ for\footnote{
i.e., $\big\Vert\mathcal S^\varphi(e,e)\big\Vert_\mathcal M\le\alpha\,\sqrt{\strut h(\sqrt\gamma)}\vert e\vert_M^2$
for all $e\in TM$}
some $\alpha>0$, outside a compact subset of $M$.
\end{itemize}

Then, $(h,\gamma)$ is an OY-pair for the Hessian in $(M,g^M)$. An analogous statement holds in the case
of OY-pairs for the Laplacian, with (c) replaced by

\begin{itemize}
\item[(c')] $\big\vert\vec H^\varphi\big\vert\le\alpha\,\sqrt{\strut h(\sqrt\gamma)}$ for
some $\alpha>0$, outside a compact subset of $M$.
\end{itemize}
\end{theorem}
\begin{proof}The function $h$ by hypothesis satisfies the conditions (1)--(3) of Theorem \ref{thm:sufcondOY} so we only  need to show that the function $\gamma$ satisfies the conditions  (4)--(6). The
assumptions (a) and (b) clearly imply that $\gamma$ is proper. For if $p_n\in M$ is a divergent sequence, $\rm{dist}_{M}(p_{n}, p_{0})\to \infty$ as $n\to \infty$, then  $\rm{dist}_{\mathcal{M}}(\varphi(p_{n}), \varphi(p_{0}))\to \infty$ and since $\pi_\mathcal X\big(\varphi(M)\big)$ is contained in a compact subset $K$ of $\mathcal X$ we have that $\rm{dist}_{M}(\Gamma^\mathrm v\circ \varphi(p_{n}), \Gamma^\mathrm v\circ \varphi (p_{0}))\to \infty$ as $n\to \infty$ since $\varphi$ and $\Gamma^\mathrm v$ are proper. This proves that $\gamma $ is proper, (condition (4)).

For $\xi\in T_{\varphi(p)}\mathcal M$, we write $\xi=\xi^\mathrm t+\xi^\perp$, where
$\xi^\mathrm t\in\mathrm{Im}\big(\mathrm d\varphi(p)\big)$ and $\xi^\perp\in\mathrm{Im}\big(\mathrm d\varphi(p)\big)^\perp$.
Given $p\in M$, we have using \eqref{eq:gradcomposition}, \eqref{eq:gradv} that  $$\big\vert\mathrm{grad}^M\gamma(p)\big\vert_M \le \frac{1}{\psi^2}\big\vert\mathrm{grad}^\mathcal V\Gamma\big(\pi_\mathcal V(\varphi(p))\big)
\big\vert_\mathcal V\le c\,\sqrt{\Gamma\big(\pi_\mathcal V(\varphi(p))\big)}= c\,\sqrt{\gamma(p)} $$
This shows condition (5) of Theorem \ref{thm:sufcondOY}.
Moreover, let $p\in M$ and $e\in T_pM$

\begin{multline*}
\mathrm{Hess}^M\gamma(e,e) \stackrel{\text{by \eqref{eq:hessianimmersionvertical-1}}}=
-2\,g^\mathcal X\big(\frac{\mathrm{grad}^\mathcal X\psi}{\psi},\mathrm e^\mathrm{hor}\big)\cdot
g^\mathcal V\big(\mathrm{grad}^\mathcal V\Gamma,\mathrm e^\mathrm{ver}\big) \\  +
\mathrm{Hess}^\mathcal V\Gamma\big(\mathrm e^\mathrm{ver},\mathrm e^\mathrm{ver}\big)+
g^\mathcal V\big(\mathrm{grad}^\mathcal V\Gamma,\mathcal S^\varphi(e,e)^\mathrm{ver}\big)
\\  \le   2A_0\big\vert\mathrm{grad}^\mathcal V\Gamma\big\vert_\mathcal V\cdot\vert e\vert_M^2+c'\sqrt{\gamma\cdot
h(\gamma)^\frac12}\cdot\big\vert\mathrm e^\mathrm{ver}\big\vert_\mathcal V^2\\ +
\big\vert\mathrm{grad}^\mathcal V\Gamma\big\vert_\mathcal V\cdot\big\vert\mathcal S^\varphi(e,e)^\mathrm{ver}\big\vert_\mathcal V
\\ \le  2A_0\big\vert\mathrm{grad}^\mathcal V\Gamma\big\vert_\mathcal V\cdot\vert e\vert_M^2+c'\sqrt{\gamma\cdot
h(\gamma)^\frac12}\cdot\frac{1}{\psi^2}\big\vert\mathrm e\big\vert_\mathcal M^2\\ +
\big\vert\mathrm{grad}^\mathcal V\Gamma\big\vert_\mathcal V\cdot\frac{1}{\psi}\big\vert\mathcal S^\varphi(e,e)\big\vert_\mathcal M
\\ \le \left[2A_0\big\vert\mathrm{grad}^\mathcal V\Gamma\big\vert_\mathcal V+\frac{c'}{B_0^2}
\sqrt{\strut\gamma\cdot h(\sqrt\gamma)}+\frac{\big\vert\mathrm{grad}^\mathcal V\Gamma\big\vert_\mathcal V}{B_0}\sqrt{\strut
h(\sqrt\gamma)}\right]\,\big\vert e\vert_M^2\\
\le \left[2A_0\sqrt\gamma+\left(\frac{c'}{B_0^2}+\frac1{B_0}\right)\,
\sqrt{\strut\gamma\cdot h(\sqrt\gamma)}\right]\,\big\vert e\vert_M^2,
\end{multline*}
where
\[A_0=\max_K\big\vert\frac{\mathrm{grad}^\mathcal X\psi}{\psi}\big\vert,\qquad B_0=\min_K\psi.\]

Since $h$ is unbounded and $\gamma$ is proper, then outside a compact subset of $M$
the inequality $\gamma\le\gamma\cdot h\big(\sqrt\gamma\big)$ holds. Hence, from the last inequality we get that
there exists a positive constant $c''$ such that, outside a compact set of $M$:
\[\mathrm{Hess}^M\gamma\le c''\,\sqrt{\strut\gamma\cdot h(\sqrt\gamma)}.\]
This proves that $(h,\gamma)$ is an OY-pair for the Hessian in $(M,g^M)$.
The statement for the Laplacian is proved similarly.
\end{proof}
\begin{remark}\label{thm:remOYLaplacian}
Observe that for the last statement of Theorem~\ref{thm:OYimmersion}, concerning the validity
of the Omori--Yau principle for the Laplacian in $(M,g^M)$, it is necessary to assume that
$(h,\Gamma)$ is an OY-pair for the Hessian in $(\mathcal V,g^\mathcal V)$. We also observe
that assumption (c') can be weakened by requiring that only the vertical component
of $\vec H^\varphi$ has norm less than or equal to $\alpha\,\sqrt{\strut h(\sqrt\gamma)}$ outside
some compact set.
\end{remark}

\begin{cor}
Under the assumptions of Theorem~\ref{thm:OYimmersion}, the Omori--Yau Maximum Principle
holds for $(M,g^M)$.
\end{cor}
\end{section}
\begin{section}{Curvature estimates}
We will generalize the results of  \cite{AlBeDa09} and \cite{AlBeMo09}
to the case of isometric immersions into warped products. For this,
let us consider an isometric immersion $\varphi\colon M\to\mathcal M=\mathcal X\times_\psi\mathcal V$
of the Riemannian manifold $(M,g^M)$ into a warped product manifold $\mathcal M $. Set $n_M=\mathrm{dim}(M)$
and suppose that  $n_M\ge n_\mathcal V+1$.
We will assume that there exists a point
$x_0\in\mathcal X$, a real number $b$ and a positive number $r<\mathrm{inj}_{\mathcal X}(x_0)$ such that the radial sectional curvatures $K_{\mathcal X_{x_0}}$  along the radial geodesics issuing from $x_0$ satisfies
\begin{equation}\label{eq:assradialcurvature}
\phantom{\quad\text{in $B_\mathcal X(x_0;r)$.}}K_{\mathcal X_{x_0}}\le b,\quad\text{in $B_\mathcal X(x_0;r)$.}
\end{equation}
Here $B_\mathcal X(x_0;r)$ is the geodesic ball in $\mathcal X$ centered at $x_0$ and of radius $r>0$ and  $\mathrm{inj}_{\mathcal X}(x_0)$ is the \emph{injectivity radius} of $\mathcal X$ at $x_0$.
Our estimates will be given in terms of the function $C_b$, defined by
\begin{equation}\label{eq:defCb}
C_b(t)=\begin{cases}\sqrt b\,\cot\big(\sqrt b\,t\big),&\text{if $b>0$ and $t\in(0,\pi/2\sqrt b)$}\\[.3cm]
\dfrac1t&\text{if $b=0$ and $t>0$}\\[.3cm]
\sqrt{-b}\coth\big(\sqrt{-b}\,t\big)&\text{if $b<0$ and $t>0$.}\end{cases}
\end{equation}
Observe that $C_b$ is strictly decreasing in its domain.
Denote by $\rho:\mathcal X\to\mathbb R$ the function
\[\rho(x)=\mathrm{dist}_\mathcal X(x_0,x);\]
this is a smooth function in $B_\mathcal X(x_0;r)$. The gradient of $\rho$ satisfies
\begin{equation}\label{eq:gradrhonorm1}
\big\vert\mathrm{grad}^\mathcal X\rho\big\vert_\mathcal X=
\big\vert\mathrm{grad}^\mathcal M\rho^\mathrm h\big\vert_\mathcal M\equiv1.
\end{equation}
By the \emph{Hessian Comparison Theorem} (see for instance
Ref.~\cite{GreWu}),
given $x\in B_\mathcal X(x_0;r)$ and a vector $X\in T_x\mathcal X$ orthogonal to $\mathrm{grad}^\mathcal X\rho(x)$,
then
\begin{equation}\label{eq:hessiancomparison}
\mathrm{Hess}^\mathcal X\rho(X,X)\ge C_b\big(\rho(x)\big)\vert X\vert_\mathcal X^2;
\end{equation}
on the other hand, if $Y\in T_x\mathcal X$ is parallel to $\mathrm{grad}^\mathcal X\rho(x)$
\begin{equation}\label{eq:hessianrhozero}
\mathrm{Hess}^\mathcal X\rho(X,Y)=\mathrm{Hess}^\mathcal X\rho(Y,Y)=0.
\end{equation}
\subsection{Sectional curvature estimates}
We will first generalize the main result in \cite{AlBeMo09} to the case of isometric immersions
into warped products. As above, let $\varphi\colon M\hookrightarrow\mathcal M=\mathcal X\times_\psi\mathcal V$
be an isometric immersion, let $x_0$ be a point in $\mathcal X$, denote by $\rho\colon \mathcal X\to\mathbb R$
the function $\rho(x)=\mathrm{dist}^\mathcal X(x_0,x)$, and define $F\colon \mathcal X\to\mathbb R$ by
\[F=\phi_b\circ\rho,\]
where $\phi_b$ is the function
\begin{equation}\label{eq:defphib}
\phi_b(t)=\begin{cases}1-\cos\big(\sqrt b\,t\big),&\text{if $b>0$ and $t\in(0,\pi/2\sqrt b)$}
\\[.3cm]
t^2,&\text{if $b=0$ and $t>0$}\\[.3cm]
\cosh\big(\sqrt{-b}\,t\big),&\text{if $b<0$ and $t>0$.}
\end{cases}
\end{equation}
This is a strictly increasing function, as $\phi_{b}'(t)>0$ for all $t$ in its domain, and it satisfies the following
ordinary differential equation
\begin{equation}\label{eq:EDOphibCb}
\phi_b''(t)-\phi_b'(t)C_b(t)=0.
\end{equation}
We assume that $x$ belongs to a sufficiently small
neighborhood $\mathcal U$ of $x_0$, so that $F$ is a smooth function, and that the image $\pi_\mathcal X\big(\varphi(M)\big)$
is contained in such neighborhood. The value of the parameter $b$ is chosen in such a way that inequality
\eqref{eq:assradialcurvature} holds in $\mathcal U$.
Thus we have a smooth function $f\colon M\to\mathbb R$ defined by
\begin{equation}\label{eq:defsecondaf}
f=F^\mathrm h\circ\varphi.
\end{equation}
Given $p\in M$, set $\varphi(p)=(x,v)\in\mathcal M$;
the gradient of $f$ at $p$ is computed easily from the formula
\begin{equation}\label{eq:luis1}
\phi_b'\big(\rho(x)\big)\,\mathrm{grad}^\mathcal X\rho(x)=\mathrm{grad}^\mathcal MF^\mathrm h\big(\varphi(p)\big)=
\mathrm{grad}^Mf(p)+\mathrm{grad}^\mathcal MF^\mathrm h\big(\varphi(p)\big)^\perp.
\end{equation}
Moreover, using \eqref{eq:hessianimmersionhorizontal}, for
$e\in T_pM$ one computes  the Hessian
\begin{eqnarray}\label{eq:luis2}
\mathrm{Hess}^Mf(e,e)& =& \mathrm{Hess}^\mathcal X(\phi_b\circ\rho)\big(\mathrm e^\mathrm{hor},\mathrm e^\mathrm{hor}\big) \nonumber \\ && \nonumber
\\&& +\, \phi_b'\big(\rho(x)\big)\,g^\mathcal M\big(\mathrm{grad}^\mathcal X\rho,\mathcal S^\varphi(e,e)\big) \\ && \nonumber \\
&& +\,
\phi_b'\big(\rho(x)\big)\,g^\mathcal X\big(\mathrm{grad}^\mathcal X\rho,\tfrac{\mathrm{grad}^\mathcal X\psi}\psi\big)\,
g^\mathcal M\big(\mathrm e^\mathrm{ver},\mathrm e^\mathrm{ver}\big).\nonumber
\end{eqnarray}
Moreover,
\begin{multline}\label{eq:luis3}
\mathrm{Hess}^\mathcal X(\phi_b\circ\rho)\big(\mathrm e^\mathrm{hor},\mathrm e^\mathrm{hor}\big)\\=
\phi_b''(\rho)\,g^\mathcal X\big(\mathrm{grad}^\mathcal X\rho,\mathrm e^\mathrm{hor}\big)^2+\,\phi'_b(\rho)\,
\mathrm{Hess}^\mathcal X\rho\big(\mathrm e^\mathrm{hor},\mathrm e^\mathrm{hor}\big)\\
\stackrel{\text{by \eqref{eq:EDOphibCb}}}=
\phi_b'(\rho)C_b(\rho)g^\mathcal X\big(\mathrm{grad}^\mathcal X\rho,\mathrm e^\mathrm{hor}\big)^2 +\,\phi'_b(\rho)\,
\mathrm{Hess}^\mathcal X\rho\big(\mathrm e^\mathrm{hor},\mathrm e^\mathrm{hor}\big)
\\=\phi_b'(\rho)\left[C_b(\rho)g^\mathcal X\big(\mathrm{grad}^\mathcal X\rho,\mathrm e^\mathrm{hor}\big)^2 +\,
\mathrm{Hess}^\mathcal X\rho\big(\mathrm e^\mathrm{hor},\mathrm e^\mathrm{hor}\big)\right].
\end{multline}
Let us now recall the following result, known in the literature as \emph{Otsuki's Lemma}.
\begin{lemma}
Let $\beta:V\times V\to W$ be a symmetric bilinear form, where $V$ and $W$ are finite dimensional
vector spaces. Assume that $\beta(v,v)\ne0$ for all $v\ne0$, and that $\mathrm{dim}(W)<\mathrm{dim}(V)$.
Then, there exist linearly independent vectors $v_1,v_2\in V$ such that $\beta(v_1,v_1)=\beta(v_2,v_2)$
and $\beta(v_1,v_2)=0$.
\end{lemma}
\begin{proof}
See for instance \cite[page 28]{KobaNomi}.
\end{proof}
We can now state our first main result (Theorem~A in the Introduction).
\begin{theorem}\label{thm:mainsectcurvestimate}
Let $(M,g^M)$ be a complete Riemannian $n_{M}$-manifold for which the Weak Omori--Yau principle for the Hessian holds and
let $\varphi:M\to\mathcal M=\mathcal X\times_\psi\mathcal V$ be an isometric immersion.
Assume that the following hypotheses are satisfied
\begin{enumerate}
\item\label{itm:sectcurv2} $\pi_\mathcal X\big(\varphi(M)\big)\subset B_\mathcal X(x_0;r)$ for some $x_0\in\mathcal X$ and $r\in\left(0,
\mathrm{inj}_\mathcal X(x_0)\right)$.
\item[]
\item Assumption \eqref{eq:assradialcurvature} holds.
\item[]
\item\label{itm:sectcurv3} $2n_M\ge 2n_\mathcal V+n_\mathcal X+1$.
\end{enumerate}
Then,
\begin{equation}\label{eq:sectcurvestimate}
\sup_MK_M\ge C_b(r)^2+\inf_{B_\mathcal X(x_0;r)}K_\mathcal X.
\end{equation}
\end{theorem}
\begin{proof}The assumption \eqref{itm:sectcurv3} together with the natural  dimension bound $n_M\leq  n_\mathcal M-1$
 implies that
\[2n_\mathcal V+n_\mathcal X+1\le2n_M\le2n_\mathcal X+2n_\mathcal V-2.\]
And that gives
\begin{equation}\label{eq:boundnX}
n_\mathcal X\ge3.
\end{equation}
Using again assumption \eqref{itm:sectcurv3} and \eqref{eq:boundnX}, together with the fact that
$\mathrm d\varphi$ is injective, we have that for all $p\in M$ there exists a subspace $\Pi_p\subset T_pM$
with dimension
\begin{equation}\label{eq:dimPi}
\mathrm{dim}(\Pi_p)\ge n_M-n_\mathcal V\ge\frac{1}{2}\big(n_\mathcal X+1\big)\ge2,
\end{equation}
such that $\mathrm d\varphi(\Pi_p)$ is horizontal. Thus, for all $e\in\Pi_p$,
$\mathrm e^\mathrm{ver}=0$, and
\[\mathrm e=\mathrm e^\mathrm{hor}=g^\mathcal X\big(\mathrm e^\mathrm{hor},
\mathrm{grad}^\mathcal X\rho(p)\big)\mathrm{grad}^\mathcal X\rho+\mathrm e^\perp.\]
Using Hessian's comparison theorem, we get
\begin{eqnarray}\label{eq:luis4}
\mathrm{Hess}^\mathcal X\rho\big(\mathrm e^\mathrm{hor}, \mathrm e^\mathrm{hor}\big)& = &
\mathrm{Hess}^\mathcal X\rho\big(\mathrm e^\perp, \mathrm e^\perp\big)\nonumber \\ && \nonumber \\
& \ge & C_b(\rho)\,
g^\mathcal X\big(\mathrm e^\perp,\mathrm e^\perp\big) \\ && \nonumber \\& = & C_b(\rho)\Big[\big\vert\mathrm e\big\vert_\mathcal M^2-g^\mathcal X
\big(\mathrm e^\mathrm{hor},\mathrm{grad}^\mathcal X\rho\big)^2\Big].\nonumber
\end{eqnarray}
From \eqref{eq:luis3}, we obtain that for all $e\in\Pi_p$
\begin{equation}\label{eq:luis5}
\mathrm{Hess}^\mathcal X(\phi_b\circ\rho)\big(\mathrm e^\mathrm{hor},\mathrm e^\mathrm{hor}\big)
\ge\phi_b'(\rho)\,C_b(\rho)\,\big\vert\mathrm e\big\vert_\mathcal M^2.
\end{equation}
Moreover, for the function $f:M\to\mathbb R$ defined in \eqref{eq:defsecondaf},
from \eqref{eq:luis2} and \eqref{eq:luis5}, we obtain
\begin{equation}\label{eq:luis6}
\mathrm{Hess}^Mf(e,e)\ge\phi_b'(\rho)\Big[C_b(\rho)\,\big\vert\mathrm e\big\vert_\mathcal M^2
-\big\vert\mathcal S^\varphi(e,e)\big\vert_\mathcal M\Big],
\end{equation}
for all $p\in M$ and all $e\in\Pi_p$.
Here, we have used the equalities $\big\vert\mathrm{grad}^\mathcal X\rho\big\vert_\mathcal X=1$
and $\mathrm e^\mathrm{ver}=0$.
We will now apply the Omori--Yau principle for the Hessian to the function $f$, which is
smooth and bounded by assumption \eqref{itm:sectcurv2}. Let $p_n\in M$ be a sequence satisfying
\begin{itemize}
\item[(a)] $f(p_n)>\sup_Mf-\displaystyle\frac{1}{n}$;\smallskip

\item[(b)] $\mathrm{Hess}^Mf(p_n)<\displaystyle\frac{1}{n}$,
\end{itemize}
for all $n$. Choose $e\in\Pi_{p_n}$, and recall that $\vert e\vert_M=\big\vert\mathrm e\big\vert_\mathcal M=
\big\vert\mathrm e^\mathrm{hor}\big\vert_\mathcal M$.
By (b) and \eqref{eq:luis6}, we have
\[\frac{1}{n}\vert e\vert_M^2>\mathrm{Hess}^Mf(p_n)(e,e)\ge\phi_b'(s_n)\left(C_b(s_n)\,\vert e\vert_M^2-\big\vert\mathcal S^\varphi(e,e)\big\vert_\mathcal M\right),\]
where
\[s_n=\rho^\mathrm h\big(\varphi(p_n)\big).\]
Hence
\begin{equation}\label{eq:luis7}
\big\vert\mathcal S^\varphi(e,e)\big\vert_\mathcal M\ge\left(C_b(s_n)-\frac{1}{n\,\phi_b'(s_n)}\right)\vert e\vert_M^2
\ge\left(C_b(r)-\frac{1}{n\,\phi_b'(s_n)}\right)\vert e\vert_M^2.
\end{equation}
We now observe that assumption \eqref{itm:sectcurv3} gives $n_M>n_\mathcal V$; this implies that the
image $\varphi(M)$ is not contained in the vertical fiber $\{x_0\}\times\mathcal V$, and in particular
that $\sup_M\rho^\mathrm h\circ\varphi>0$. For all $b$, the function $\phi_b$ is increasing and positive,
and therefore $\sup_Mf>0$; this says that $\pi_\mathcal X(p_n)$ stays away from $x_0$, i.e., there exists
$\delta>0$ such  that $s_n\ge\delta$. Therefore
\begin{equation}\label{eq:limphib'}
\lim_{n\to\infty}\frac1{n\,\phi_b'(s_n)}=0,
\end{equation}
and so for $n$ sufficiently large, we have $C_b(s_n)-\displaystyle\frac{1}{n\,\phi_b'(s_n)}>0$, which implies in particular
that $\mathcal S^\varphi(e,e)\ne0$ for all $e\in\Pi_{p_n}\setminus\{0\}$. We can invoke Otsuki's Lemma, applied to
the symmetric bilinear form given by the restriction of $\mathcal S^\varphi$ to $\Pi_{p_n}\times\Pi_{p_n}$,
which takes values in the space $\mathrm{Im}\big(\mathrm d\varphi(p_n)\big)^\perp$.
Note that, by assumption \eqref{itm:sectcurv3}
\[\mathrm{dim}\left[\mathrm{Im}\big(\mathrm d\varphi(p_n)\big)^\perp\right]=n_\mathcal X+n_\mathcal V-n_M<
n_M-n_\mathcal V=\mathrm{dim}(\Pi_{p_n}).\]
Thus, there exist linearly independent vectors $e^1,e^2\in\Pi_{p_n}$ such that:
\[\mathcal S^\varphi(e^1,e^1)=\mathcal S^\varphi(e^2,e^2),\quad\mathcal S^\varphi(e^1,e^2)=0.\]
We can assume $\vert e^1\vert_M\ge\vert e^2\vert_M>0$.
We will now compare the sectional curvature $K_M(e^1,e^2)$  with the sectional curvature $K_\mathcal M\big(\mathrm d\varphi(e^1),\mathrm d\varphi(e^2)\big)$.
The  plane spanned by $ \mathrm d\varphi(e^1),\,\mathrm d\varphi(e^2)$ is horizontal, and recalling \eqref{eq:equalitysectcurvature}, we have:
\begin{equation}\label{eq:equalitysectcurvatureswarped}
K_\mathcal M\big(\mathrm d\varphi(e^1),\mathrm d\varphi(e^2)\big)=
K_\mathcal X\big(\mathrm d\varphi(e^1),\mathrm d\varphi(e^2)\big).
\end{equation}
Then, using Gauss equation we have
\begin{eqnarray}
K_M(e^1,e^2)-K_\mathcal X\big(\mathrm d\varphi(e^1),\mathrm d\varphi(e^2)\big)&=&
K_M(e^1,e^2)-K_\mathcal M\big(\mathrm d\varphi(e^1),\mathrm d\varphi(e^2)\big) \nonumber\end{eqnarray}
\begin{eqnarray}
= \frac{g^\mathcal M\big(\mathcal S(e^1,e^1),\mathcal S(e^2,e^2)\big)-\big\vert\mathcal S(e^1,e^2)\big\vert_\mathcal M^2}{%
\vert e^1\vert^2_M\vert e^2\vert_M^2-g^M(e^1,e^2)^2}& = &\frac{\big\vert\mathcal S(e^1,e^1)\big\vert_\mathcal M^2}{%
\vert e^1\vert^2_M\vert e^2\vert_M^2-g^M(e^1,e^2)^2}\nonumber \\ \geq  \left[\frac{\big\vert\mathcal S(e^1,e^1)\big\vert_\mathcal M}{
\vert e^1\vert^2_M}\right]^2 & \stackrel{\text{by \eqref{eq:luis7}}}{\geq} &  \left(C_b(r)-\frac{1}{n\,\phi_b'(s_n)}\right)^2.\nonumber
\end{eqnarray}
Hence
\[\sup_MK_M-\inf_{B_\mathcal X(x_0;r)}K_\mathcal X\ge\left(C_b(r)-\frac{1}{n\,\phi_b'(s_n)}\right)^2.\]
Taking the limit as $n\to\infty$ and recalling \eqref{eq:limphib'} we get
\eqref{eq:sectcurvestimate}.
\end{proof}

\subsection{Mean curvature estimates}
We start with an elementary result.
\begin{lemma}\label{thm:lemmaHS}
Let $\big(W_i,\langle\,\,,\,\rangle_i\big)$, $i=1,2$, be finite dimensional vector spaces
with inner product, with dimensions $n_i$, $i=1,2$, and let $T:W_1\to W_2$ be a linear map with the property
that there exists an orthogonal decomposition $W_1=W\oplus W'$ such that $T\vert_W:W\to W_2$ is a surjective
isometry and $T\vert_{W'}=0$. Then, for every orthonormal basis $\xi_1,\ldots,\xi_{n_1}$ of $W_1$, the
following equality holds.
\begin{equation}\label{eq:sumequaln2}
\sum_{i=1}^{n_1}\Big\vert T\xi_i\big\vert_2^2=n_2.
\end{equation}
In particular, if $\eta_1,\ldots,\eta_n$ is any orthonormal family in $W_1$, then $\sum\limits_{i=1}^n\big\vert T\eta_i\big\vert_2^2\le n_2$.
\end{lemma}
\begin{proof}
The left-hand side of \eqref{eq:sumequaln2} does not depend\footnote{%
Namely, if $\eta_1',\ldots,\eta_{n_1}'$ is another orthonormal basis, then there exists an orthogonal
$n_1\times n_1$ matrix $A=(a_{ij})$ such that $\eta_i'=\sum_ja_{ij}\eta_j$ for all $i$. The orthogonality of $A$ means
that $(A^*A)_{jk}=\sum_ia_{ij}a_{ik}=\delta_{jk}$ for all $j,k$. Then:
\[\sum_i\vert T\eta'_i\vert_2^2=\sum_{i,j,k}a_{ij}a_{ik}\langle T\eta_j,T\eta_k\rangle_2=
\sum_{j,k}\delta_{jk}\langle T\eta_j,T\eta_k\rangle_2=\sum_j\vert T\eta_j\vert_2^2.\]}
on the orthonormal basis; it is the \emph{Hilbert--Schmidt}
squared norm of the linear map $T$. The equality is verified easily using an orthonormal basis of $W_1$ consisting
of the union of an orthonormal basis of $W$ and an orthonormal basis of $W'$.
\end{proof}
We can now prove the following (Theorem~B in the Introduction):
\begin{theorem}\label{thm:meancurvimmwarpedpr}
Let $(M,g^M)$ be a complete  Riemannian manifold for which the Weak Omori--Yau principle for the Laplacian holds, and
let $\varphi:M\to\mathcal M=\mathcal X\times_\psi\mathcal V$ be an isometric immersion.
Assume that the following hypotheses are satisfied.
\smallskip
\begin{enumerate}
\item\label{itm:sectcurv1} $\pi_\mathcal X\big(\varphi(M)\big)\subset B_\mathcal X(x_0;r)$ for some $x_0\in\mathcal X$ and $r\in(0,
\mathrm{inj}_\mathcal X(x_0))$.
\item[]
\item Assumption \eqref{eq:assradialcurvature} holds.
\end{enumerate}
Then, denoting by $\vec H^\varphi$ the mean curvature vector of $\varphi$, one has the following estimate
on the supremum of $\vert\vec H^\varphi\big\vert_\mathcal M$.
\begin{equation}\label{eq:meancurvimmwarpedpr}
\sup_M\big\vert\vec H^\varphi\big\vert_\mathcal M\ge(n_M-n_\mathcal V)\, C_b(r)-{n_\mathcal V}\Psi_0,
\end{equation}
where
\begin{equation}\label{eq:defPsi0}
\Psi_0=\sup_{\mathrm{dist}(x,x_0)\le r}\Big\vert\frac{\mathrm{grad}^\mathcal X\psi(x)}{\psi(x)}\Big\vert_\mathcal X.
\end{equation}
\end{theorem}
\begin{proof}
Inequality \eqref{eq:meancurvimmwarpedpr} is proved applying the Weak Omori--Yau principle to the
function $f:M\to\mathbb R$ defined in \eqref{eq:defsecondaf}. Let us give an estimate for the Laplacian
of $f$ as follows.
From \eqref{eq:luis2} and \eqref{eq:luis3}, given $p\in M$ and $e\in T_pM$ we have
\begin{eqnarray}\label{eq:new1}
\mathrm{Hess}^Mf(e,e) \!\!\!&=&\!\! \!\phi_b'(\rho)\left[C_b(\rho) g^\mathcal X\big(\mathrm{grad}^\mathcal X\rho,\mathrm e^\mathrm{hor}\big)^2
+\mathrm{Hess}^\mathcal X\rho\big(\mathrm e^\mathrm{hor},\mathrm e^\mathrm{hor}\big)\right]\nonumber
\\&& \!\! +\,
\phi_b'(\rho) \,g^\mathcal X\big(\mathrm{grad}^\mathcal X\rho,\frac{\mathrm{grad}^\mathcal X\psi}\psi\big)
\,g^\mathcal M\big(\mathrm e^\mathrm{ver},\mathrm e^\mathrm{ver}\big) \\ &&
\!\!+\,\phi_b'(\rho)\, g^\mathcal M\big(\mathrm{grad}^\mathcal X\rho,\mathcal S^\varphi(e,e)\big).\nonumber
\end{eqnarray}
Let us write
$\mathrm e=\mathrm e^\mathrm{hor}+\mathrm e^\mathrm{ver}$
and
$ e^\mathrm{hor}=\mathrm e^\rho+\mathrm e^\perp,$
where \[\mathrm e^\rho=g^\mathcal X\big(\mathrm e^\mathrm{hor},\mathrm{grad}^\mathcal X\rho\big)\,
\mathrm{grad}^\mathcal X\rho.\]

Observe
$$
\big\vert\mathrm e^\mathrm{hor}\big\vert_\mathcal M^2=
\big\vert\mathrm e^\rho\big\vert_\mathcal M^2+\big\vert\mathrm e^\perp\big\vert_\mathcal M^2=
g^\mathcal X\big(\mathrm e^\mathrm{hor},\mathrm{grad}^\mathcal X\rho\big)^2+\big\vert\mathrm e^\perp
\big\vert_\mathcal M^2
$$
Using the Hessian Comparison Theorem \eqref{eq:hessiancomparison}, we obtain
\begin{eqnarray}\label{eq:new2}
\mathrm{Hess}^\mathcal X\rho\big(\mathrm e^\mathrm{hor},\mathrm e^\mathrm{hor}\big)&=&
\mathrm{Hess}^\mathcal X\rho\big(\mathrm e^\perp,\mathrm e^\perp\big)\nonumber \\
& \ge &  C_b(\rho)\,g^\mathcal X\big(\mathrm e^\perp,\mathrm e^\perp\big) \\ & = &
C_b(\rho)\big\vert\mathrm e^\perp\big\vert_\mathcal M^2.\nonumber
\end{eqnarray}
From \eqref{eq:new1} and \eqref{eq:new2} we get to the following inequality.
\begin{eqnarray}\label{eq:new3}
\mathrm{Hess}^Mf(e,e)& \ge & \phi_b'(\rho)\,C_b(\rho)\left[g^\mathcal X\big(\mathrm e^\mathrm{hor},\mathrm{grad}^\mathcal X\rho\big)^2+\big\vert\mathrm e^\perp
\big\vert_\mathcal M^2\right]\nonumber \\
 && \nonumber \\
 && +
\phi_b'(\rho) g^\mathcal M\big(\mathrm{grad}^\mathcal X\rho,\frac{\mathrm{grad}^\mathcal X\psi}\psi\big)
\,\big\vert\mathrm e^\mathrm{ver}\big\vert_\mathcal M^2\nonumber \\
&& \nonumber \\
&& +\phi_b'(\rho) g^\mathcal M\big(\mathrm{grad}^\mathcal X\rho,\mathcal S^\varphi(e,e)\big)\nonumber \\
&&\nonumber\\
& = &\phi_b'(\rho)\,C_b(\rho)\big\vert\mathrm e^\mathrm{hor}
\big\vert_\mathcal M^2 \\
&&\nonumber\\
&& +
\phi_b'(\rho) g^\mathcal M\big(\mathrm{grad}^\mathcal X\rho,\frac{\mathrm{grad}^\mathcal X\psi}\psi\big)
\,\big\vert\mathrm e^\mathrm{ver}\big\vert_\mathcal M^2\nonumber \\
&&\nonumber \\
&& +\phi_b'(\rho) g^\mathcal M\big(\mathrm{grad}^\mathcal X\rho,\mathcal S^\varphi(e,e)\big)\nonumber
\end{eqnarray}

\smallskip
Let $(e_i)_{i=1}^{n_M}$ be an orthonormal basis of $T_pM$; from \eqref{eq:new3}, we get
\begin{eqnarray}\label{eq:new4}
\triangle^Mf& = &\sum_{i=1}^{n_M}\mathrm{Hess}^Mf(e_i,e_i)\nonumber \\
&&\nonumber \\
&\ge & \phi_b'(\rho)\,C_b(\rho)\sum_{i=1}^{n_M}\big\vert
\mathrm d\varphi(e_i)^\mathrm{hor}
\big\vert_\mathcal M^2 \nonumber \\
&&\nonumber \\
&& +
\phi_b'(\rho) g^\mathcal M\big(\mathrm{grad}^\mathcal X\rho,\frac{\mathrm{grad}^\mathcal X\psi}\psi\big)
\sum_{i=1}^{n_M}\big\vert\mathrm d\varphi(e_i)^\mathrm{ver}\big\vert_\mathcal M^2\nonumber \\
&&\nonumber \\
&& +\phi_b'(\rho) g^\mathcal M\big(\mathrm{grad}^\mathcal X\rho,\vec H^\varphi\big)\nonumber \\
&&\nonumber \\
&\ge& \phi_b'(\rho)\,C_b(\rho)\sum_{i=1}^{n_M}\left(1-\big\vert
\mathrm d\varphi(e_i)^\mathrm{ver}
\big\vert_\mathcal M^2\right) \\
&&\nonumber \\
&& -
\phi_b'(\rho) \big\vert\frac{\mathrm{grad}^\mathcal X\psi}\psi\big\vert_\mathcal X
\sum_{i=1}^{n_M}\big\vert\mathrm d\varphi(e_i)^\mathrm{ver}\big\vert_\mathcal M^2
-\phi_b'(\rho)\big\vert\vec H^\varphi\big\vert_\mathcal M \nonumber \\
&&\nonumber\\
& =& \phi_b'(\rho)\Big[C_b(\rho)\Big(n_M-\sum_{i=1}^{n_M}\big\vert
\mathrm d\varphi(e_i)^\mathrm{ver}
\big\vert_\mathcal M^2\Big)\nonumber \\
&&\nonumber \\
&&-\big\vert\frac{\mathrm{grad}^\mathcal X\psi}\psi\big\vert_\mathcal X
\sum_{i=1}^{n_M}\big\vert\mathrm d\varphi(e_i)^\mathrm{ver}\big\vert_\mathcal M^2
-\big\vert\vec H^\varphi\big\vert_\mathcal M\Big].\nonumber
\end{eqnarray}

Now, we claim that the following inequality holds
\begin{equation}\label{eq:HilbSchm}
\sum_{i=1}^{n_M}\big\vert e_i^\mathrm{ver}\big\vert^2_\mathcal M\le n_\mathcal V.
\end{equation}
This follows from Lemma~\ref{thm:lemmaHS} applied to the linear map
\[\mathrm d\pi_\mathcal V\big(\varphi(p_n)\big):T_{\varphi(p_n)}\mathcal M\longrightarrow T_{\pi(\varphi(p_n))}\mathcal V,\]
where the space $T_{\pi(\varphi(p_n))}\mathcal V$ is endowed with the inner product $\psi^2\cdot g^\mathcal V$,
considering the orthonormal family $e_1,\ldots,e_{n_M}$ in $T_{\varphi(p_n)}\mathcal M$.

Thus, \eqref{eq:new4} gives
\begin{equation}\label{eq:new5}
\triangle^Mf\ge\phi_b'(\rho)\Big[C_b(\rho)\big(n_M-n_\mathcal V\big)-n_\mathcal V\,\Psi_0
-\big\vert\vec H^\varphi\big\vert_\mathcal M\Big].
\end{equation}
The Weak Omori--Yau Principle for the Laplacian yields the existence of
a sequence $p_n$  in $M$ such that
\begin{itemize}
\item[(a)] $f(p_n)>\sup_Mf-\displaystyle\frac{1}{n}$.
\item[]
\item[(b)] $\Delta^Mf(p_n)<\displaystyle\frac{1}{n}$.
\end{itemize}
Set $s_n=\rho\big(\varphi(p_n)\big)$.  The inequality \eqref{eq:new5} gives
\begin{equation}\label{eq:new6}
\phi_b'(s_n)\Big[C_b(s_n)\big(n_M-n_\mathcal V\big)-n_\mathcal V\,\Psi_0
-\big\vert\vec H^\varphi\big\vert_\mathcal M\Big]<\displaystyle\frac{1}{n}
\end{equation}
for all $n$.
Arguing as in the proof of Theorem~\ref{thm:mainsectcurvestimate}, the sequence $s_n$ is bounded away from
$0$, and so is the quantity $\phi_b'(s_n)$. Moreover, since $C_b$ is decreasing, it is $C_b(s_n)\ge C_b(r)$
for all $n$. Taking the limit as $n\to\infty$ in \eqref{eq:new6}, we obtain
\[C_b(r)\big(n_M-n_\mathcal V\big)-n_\mathcal V\,\Psi_0
-\sup\big\vert\vec H^\varphi\big\vert_\mathcal M\le0,\]
which is our thesis.
\end{proof}
In Theorem~\ref{thm:meancurvimmwarpedpr}, the hypothesis on the validity of the weak Omori--Yau principle
in $(M,g^M)$ can be omitted by assuming instead that the fiber of the warped product has an OY-pair and
that $\varphi$ is proper.
\begin{cor}\label{thm:cormeansectestimates}
Assume that
\begin{itemize}
\item $(\mathcal V,g^\mathcal V)$ has an OY-pair for the Hessian,
\item[]
\item $\varphi$ is proper,
\item[]
\item the ball $B_\mathcal X(x_0;r)$ has compact closure in $\mathcal X$ (for instance, if $(\mathcal X,g^\mathcal X)$
is complete).
\end{itemize}
Then the conclusion of Theorem~\ref{thm:meancurvimmwarpedpr} holds.
\end{cor}
\begin{proof}
We argue by contradiction. If \eqref{eq:meancurvimmwarpedpr} does not hold, then the norm of the mean
curvature vector $\big\vert\vec H^\varphi\big\vert_\mathcal M$
is bounded, and we can apply Proposition~\ref{thm:OYimmersion} to deduce that the (strong)
Omori--Yau principle for the Laplacian holds in $(M,g^M)$. Thus, \emph{a fortiori}, \eqref{eq:meancurvimmwarpedpr}
holds.
\end{proof}
\end{section}

\begin{section}{Isometric immersions into Riemannian submersions}
\label{sec:immersionssubmersions}
Our curvature estimates for isometric immersions into warped product can be partially extended
to the far more general case of immersions into the total space of  Riemannian submersions.
Given Riemannian smooth manifolds $(\mathcal M,g^\mathcal M)$ and $(\mathcal X,g^\mathcal X)$,
a Riemannian submersion is a smooth surjective map $\pi:\mathcal M\to\mathcal X$ such that the differential
$\mathrm d\pi$ has everywhere maximal rank, and it is an isometry when restricted to horizontal
vectors, i.e., $\vert X\vert_\mathcal M=\big\vert\mathrm d\pi(X)\vert_\mathcal X$ for all $X$ orthogonal
to the kernel of $\mathrm d\pi$. The manifold $\mathcal M$ is called the \emph{total space} of the submersion,
$\mathcal X$ is the base, and for all $x\in\mathcal X$, the \emph{fiber} $\mathcal V_x$ is the smooth embedded
submanifold of $\mathcal M$ given by $\pi^{-1}(x)$.

A horizontal vector field $X\in\mathfrak X(\mathcal M)$ is \emph{basic} if it is $\pi$-related to
some vector field $X_*\in\mathfrak X(\mathcal X)$. If $X$ and $Y$ are basic vector fields, then the
horizontal component $(\nabla^\mathcal M_XY)^\mathrm h$ of the covariant derivative $\nabla^\mathcal M_XY$
is basic, and it is $\pi$-related to $\nabla^\mathcal X_{X_*}Y_*$, see \cite[Lemma~1]{One66}.
Having this in mind, it is easy to prove the following
\begin{lemma}\label{thm:hesshorriemsubmersion}
Let $F:\mathcal X\to\mathbb R$ be a smooth function; set $F^\mathrm h=F\circ\pi\colon\mathcal M\to\mathbb R$.
Then, the gradient $\mathrm{grad}^\mathcal MF^\mathrm h$ is basic. Given $p\in\mathcal M$ and horizontal vectors
$X,Y\in T_p\mathcal M$, then $\mathrm{Hess}^\mathcal MF^\mathrm h(X,Y)$ is equal to
$\mathrm{Hess}^\mathcal XF\big(\mathrm d\pi_p(X),\mathrm d\pi_p(Y)\big)$.\qed
\end{lemma}
Let us recall that the geometry of a Riemannian submersion $\pi:\mathcal M\to\mathcal X$ is described
by the \emph{fundamental tensors} $T$ and $A$, introduced by O'Neill see \cite{One66}, defined by the following
formulas:
\[T_\xi(\eta)=\big(\nabla^\mathcal M_{\xi^\mathrm{ver}}(\eta^\mathrm{ver})\big)^\mathrm{hor}+
\big(\nabla^\mathcal M_{\xi^\mathrm{ver}}(\eta^\mathrm{hor})\big)^\mathrm{ver},\]
\[A_\xi(\eta)=\big(\nabla^\mathcal M_{\xi^\mathrm{hor}}(\eta^\mathrm{hor})\big)^\mathrm{ver}+
\big(\nabla^\mathcal M_{\xi^\mathrm{hor}}(\eta^\mathrm{ver})\big)^\mathrm{hor},\]
for $\xi,\eta\in T\mathcal M$. Restricted to vertical vectors, $T$ is the second fundamental
form of the fibers of the submersion. On horizontal fields, $A$ is essentially the integrability
tensor of the horizontal distribution of the submersion.
\subsection{Sectional curvature estimates}
The tensor $A$ is related to the sectional curvature of horizontal planes, as follows.
Let $p\in\mathcal M$ and $X,Y\in T_p\mathcal M$ be (linearly independent) horizontal vectors;
set $X_*=\mathrm d\pi_p(X)$ and $Y_*=\mathrm d\pi_p(Y)$. Then:
\begin{equation}\label{eq:sectcurvhorsubmersion}
K_\mathcal M(X,Y)=K_\mathcal X(X_*,Y_*)-\frac{3\,\big\vert A_XY\big\vert_{\mathcal M}^2}{%
\big\vert X\big\Vert_\mathcal M^2\,\big\vert Y\big\vert_\mathcal M^2-g^\mathcal M(X,Y)^2}.
\end{equation}
Given $p\in\mathcal M$, let us introduce the following notation:
\[\mathrm{sec}_\mathrm{hor}^\mathcal M(p)=\min\Big\{K_\mathcal M(\Pi):\Pi\subset T_p\mathcal M\
\text{horizontal $2$-plane}\Big\}.\]
When the horizontal distribution is integrable, i.e., when the tensor $A$ vanishes identically on
horizontal vectors, then by \eqref{eq:sectcurvhorsubmersion} the sectional curvature of horizontal
planes in $\mathcal M$ coincides with the curvature of the corresponding plane in the base manifold $\mathcal X$.
\begin{theorem}\label{thm:mainsectcurvestimatesub}
Let $(M,g^M)$ be a Riemannian manifold in which the Omori--Yau Principle for the Hessian holds, and
let $\varphi:M\to\mathcal M$ be an isometric immersion of $M$ into the total space of a Riemannian submersion
$\pi:\mathcal M\to\mathcal X$. Denote by $n_\mathcal V$ the dimension of the fibers of $\pi$.
Assume that the following hypotheses are satisfied.
\begin{enumerate}
\item\label{itm:sectcurv2sub} $\pi\big(\varphi(M)\big)\subset B_\mathcal X(x_0;r)$ for some $x_0\in\mathcal X$ and $r\in\left(0,
\mathrm{inj}_\mathcal X(x_0)\right)$
\item[]
\item assumption \eqref{eq:assradialcurvature} holds,
\item[]
\item\label{itm:sectcurv3sub} $2n_M\ge 2n_\mathcal V+n_\mathcal X+1$.
\end{enumerate}
Then,
\begin{equation}\label{eq:sectcurvestimatesub}
\sup_MK_M\ge C_b(r)^2+\inf_{\pi^{-1}(B_\mathcal X(x_0;r))}\mathrm{sec}_\mathrm{hor}^\mathcal M.
\end{equation}
If the horizontal distribution of $\pi$ is integrable, then
\begin{equation}\label{eq:sectcurvestimatesubint}
\sup_MK_M\ge C_b(r)^2+\inf_{B_\mathcal X(x_0;r)}K_\mathcal X.
\end{equation}
\end{theorem}
\begin{proof}
In the proof of Theorem~\ref{thm:mainsectcurvestimate}, the Omori--Yau
principle is used for evaluating the Hessian of smooth functions on the base $\mathcal X$, in the
directions of horizontal vectors. By Lemma~\ref{thm:hesshorriemsubmersion}, also in the case
of Riemannian submersions the value of the Hessian
in horizontal directions coincides with the Hessian on the base manifolds.
Thus, the proof of Theorem~\ref{thm:mainsectcurvestimate} can be repeated \emph{verbatim}, with the
exception of formula \eqref{eq:equalitysectcurvatureswarped}, which is replaced by \eqref{eq:sectcurvhorsubmersion}.
This yields the estimate \eqref{eq:sectcurvestimatesub}. When the horizontal distribution is integrable,
then also formula \eqref{eq:equalitysectcurvatureswarped} holds for the Riemannian submersion
$\pi:\mathcal M\to\mathcal X$, and the conclusion is exactly the same as in Theorem~\ref{thm:mainsectcurvestimate}.
\end{proof}
\subsection{Mean curvature estimates}
In order to extend to Riemannian submersions the result of Theorem~\ref{thm:meancurvimmwarpedpr},
we need a generalization of formula \eqref{eq:hessianimmersionhorizontal}.
This is obtained easily using the expressions for the Levi--Civita connection of the total
space of a Riemannian submersions given, for instance, in \cite[Lemma~3]{One66}.
If $X$ is basic, and $V$ is vertical, then the horizontal component of the covariant derivative $\nabla^\mathcal M_VX$
is given by
\begin{equation}\label{eq:nablavbasichorsub}
\big(\nabla^\mathcal M_VX\big)^\mathrm h=\big(\nabla^\mathcal M_XV\big)^\mathrm h=A_XV;
\end{equation}
similarly, the vertical component of $\nabla^\mathcal M_VX$ is
\begin{equation}\label{eq:nablavbasicversub}
\big(\nabla^\mathcal M_VX\big)^\mathrm v=T_VX.
\end{equation}
Let $F:\mathcal X\to\mathbb R$ be a smooth map, $F^\mathrm h=F\circ\pi:\mathcal M\to\mathbb R$,
$\varphi:M\to\mathcal M$ an isometric immersion and $f=F^\mathrm h\circ\varphi:M\to\mathbb R$.
For $p\in M$ and $e\in T_pM$, let us set $\xi=\mathrm d\varphi_p(e)\in T_{\varphi(p)}\mathcal M$,
and let $\mathcal S^\varphi$ denote the second fundamental form of $\varphi$..
Using the fact that $\mathrm dF^\mathrm h$ is basic, $\pi$-related to $\mathrm{grad}^\mathcal XF$,
and using formulas \eqref{eq:hessiancomposition},
\eqref{eq:nablavbasichorsub} and \eqref{eq:nablavbasicversub}, we obtain the following expression
for the Hessian of $f$.
\begin{multline}\label{eq:hessianfsubmersions}
\mathrm{Hess}^Mf(e,e)= \mathrm{Hess}^\mathcal MF^\mathrm h(\xi,\xi)+g^\mathcal M\big(\mathrm{grad}^\mathcal MF^\mathrm h,
\mathcal S^\varphi(e,e)\big) \\ =\mathrm{Hess}^\mathcal MF^\mathrm h(\xi^\mathrm{hor},\xi^\mathrm{hor})+
2\,\mathrm{Hess}^\mathcal MF^\mathrm h(\xi^\mathrm{hor},\xi^\mathrm{ver}) \\ +\mathrm{Hess}^\mathcal MF^\mathrm h(\xi^\mathrm{ver},\xi^\mathrm{ver})
+g^\mathcal M\big(\mathrm{grad}^\mathcal MF^\mathrm h,
\mathcal S^\varphi(e,e)\big)\\ = \mathrm{Hess}^\mathcal XF(\xi^\mathrm{hor}_*,\xi^\mathrm{hor}_*)+2\,
g^\mathcal M\big(A_{\xi^\mathrm{hor}}(\mathrm{grad}^\mathcal MF^\mathrm h),\xi^\mathrm{ver}\big)\\ +
g^\mathcal M\big(T_{\xi^\mathrm{ver}}(\mathrm{grad}^\mathcal MF^\mathrm h),\xi^\mathrm{ver}\big)
 +g^\mathcal M\big(\mathrm{grad}^\mathcal MF^\mathrm h,
\mathcal S^\varphi(e,e)\big).
\end{multline}
Denote by $\mathcal S^\mathcal V$ the second fundamental
form of the fibers; the term
containing the fundamental tensor $T$ can be rewritten
in terms of $\mathcal S^\mathcal V$ as
\begin{equation}\label{eq:Tintermsofsecfundform}
g^\mathcal M\big(T_{\xi^\mathrm{ver}}(\mathrm{grad}^\mathcal MF^\mathrm h),\xi^\mathrm{ver}\big)=
-g^\mathcal M\big(\mathcal S^\mathcal V(\xi^\mathrm{ver},\xi^\mathrm{ver}),\mathrm{grad}^\mathcal M F^\mathrm h\big).
\end{equation}
Comparing \eqref{eq:hessianfsubmersions} with \eqref{eq:hessianimmersionhorizontal}, the reader will observe
that the term in \eqref{eq:Tintermsofsecfundform}
corresponds to the last term in \eqref{eq:hessianimmersionhorizontal}, see \eqref{eq:secfundformfibers}.
The new term here is the one containing the fundamental tensor $A$, which vanishes in the case of warped
products. Thus, it is easy to formulate the following extension of Theorem~\ref{thm:meancurvimmwarpedpr}:
\begin{theorem}\label{thm:meancurvimmsub}
Let $(M,g^M)$ be a Riemannian manifold in which the Omori--Yau principle for the Laplacian holds, and
let $\varphi:M\to\mathcal M$ be an isometric immersion into the total space of a Riemannian submersion
$\pi:\mathcal M\to\mathcal X$.
Assume that the following hypotheses are satisfied.
\begin{enumerate}
\item $\pi\big(\varphi(M)\big)\subset B_\mathcal X(x_0;r)$ for some $x_0\in\mathcal X$ and $r\in\left]0,
\mathrm{inj}_\mathcal X(x_0)\right[$
\item assumption \eqref{eq:assradialcurvature} holds;
\end{enumerate}
Then, denoting by $\vec H^\varphi$ the mean curvature vector of $\varphi$, and by $T$, $A$ the fundamental tensors
of the Riemannian submersion, one has the following estimate
on the supremum of $\Vert\vec H^\varphi\big\Vert_\mathcal M$:
\begin{equation}\label{eq:meancurvimmsub}
\sup_M\big\Vert\vec H^\varphi\big\Vert_\mathcal M\ge(n_M-n_\mathcal V)\, C_b(r)-n_M\,\alpha_0-{n_\mathcal V}\,\tau_0,
\end{equation}
where
\begin{equation}\label{eq:deftau0}
\tau_0=\sup_{\pi^{-1}(B_\mathcal X(x_0;r))}\big\vert T\big\vert,
\end{equation}
and
\begin{equation}\label{eq:defalpha0}
\alpha_0=\sup_{\pi^{-1}(B_\mathcal X(x_0;r))}\big\vert A\big\vert.
\end{equation}
\end{theorem}
\begin{proof}
It suffices to repeat the proof of Theorem~\ref{thm:meancurvimmwarpedpr}, keeping into consideration
also the contribution of the term
\[2\,
g^\mathcal M\big(A_{\xi^\mathrm{hor}}(\mathrm{grad}^\mathcal MF^\mathrm h),\xi^\mathrm{ver}\big),\]
which is estimated as follows
\[2\,
g^\mathcal M\big(A_{\xi^\mathrm{hor}}(\mathrm{grad}^\mathcal MF^\mathrm h),\xi^\mathrm{ver}\big)\ge
-2\,\alpha_0\,\vert\xi^\mathrm{hor}\vert_\mathcal M\,\vert\xi^\mathrm{ver}\vert_\mathcal M.\]
When $\vert\xi\vert_\mathcal M^2=1$, then $\vert\xi^\mathrm{hor}\vert_\mathcal M\,\vert\xi^\mathrm{ver}\vert_\mathcal M\le\frac12$.
The conclusion follows readily.
\end{proof}

\end{section}

\end{document}